\documentclass[11pt]{amsart}

\usepackage[pdftex]{graphicx}
\usepackage{amsmath,amsfonts,amsthm,amssymb,amscd}
\usepackage[all]{xy}

\title[Hilb-desingularizations]{Hilbert desingularizations for three dimensional canonical cyclic quotient singularities}

\author[K.Sato]{Kohei Sato}
\address{National Institute Of Technology, Oyama Collage, 771 Nakakuki, Oyama, Tochigi, 323-0806, Japan}
\email{k-sato@oyama-ct.ac.jp}

\author[Y.Sato]{Yusuke Sato}
\address{Graduate School Of Mathematical Sciences, University Of Tokyo 3-8-1 Komaba, Meguro-ku, Tokyo 153-8914, Japan.}
\email{yusuke.sato@ipmu.jp}


\newtheorem{defi}{Definition}[section]
\newtheorem{thm}[defi]{Theorem}
\newtheorem{prop}[defi]{Proposition} 
\newtheorem{lem}[defi]{Lemma}

\newtheorem{cor}[defi]{Corollary}

\newtheorem{note}[defi]{Note}
\newtheorem{ex}[defi]{Example}

\renewcommand{\proofname}{Proof.\ }

\newcommand{\slmc}[1]{SL(#1,\mathbb{C})}
\newcommand{\glmc}[1]{GL(#1,\mathbb{C})}
\newcommand{\compl}[1]{\mathbb{C}^{#1}}
\newcommand{\real}[1]{\mathbb{R}^{#1}}
\newcommand{\ratio}[1]{\mathbb{Q}^{#1}}
\newcommand{\inte}[1]{\mathbb{Z}^{#1}}
\newcommand{\natu}[1]{\mathbb{N}^{#1}}

\newcommand{\GHilb}[1]{G\mathrm{\mathchar`-Hilb}(\mathbb{C}^{#1})}
\newcommand{\Hilb}{\mathrm{Hilb}}
\newcommand{\HHilb}{\mathbf{Hilb}}

\def\qed{\hfill $\Box$} 

\begin{document}

\maketitle

\begin{abstract}
In this paper, we shall discuss Hilbert property of $\Hilb^{G}(\compl{3})$, Fujiki-Oka resolutions and iterated Fujiki-Oka resolutions for three dimensional canonical cyclic quotient singularities by using the classification shown by Ishida and Iwashita\cite{II}. Finally, we shall prove that there exists a Hilbert desingularization for any three dimensional canonical cyclic quotient singularity.
\end{abstract}

\section{Introduction}\label{Intro}
Ishida and Iwashita classified three dimensional canonical cyclic quotient singularities\cite{II}. Let $G$ be a finite cyclic group of $\glmc{3}$. Then $\compl{3}/G$ has a canonical singularity if
\begin{itemize}
\item[(i)] $G \subset \slmc{3}$, 
\item[(ii)] $G=\frac{1}{r}(a_1,a_2,a_3)$ with ${\rm gcd}(a_i, a_j)=1$ and $a_i + a_j = r$ for some $1\leq i < j \leq 3$,
\item[(iii)] $G=\frac{1}{4k}(1,2k+1,4k-2)$ with $k\geq2$,
\item[(iv)]  $G=\frac{1}{9}(1,4,7)$ or $\frac{1}{14}(1,9,11)$.
\end{itemize}

The class (ii) has the two cases 
$$(\alpha)\ \frac{1}{r}(1,a,r-a)\ \text{with}\ {\rm gcd}(r,a) = 1\ \text{and}\ r>a,$$
$$(\beta)\ \frac{1}{r}(1,r-1,a)\ \text{with}\ {\rm gcd}(r,a) > 1\ \text{and}\ r>a.$$

The form (ii)-($\alpha$) is a symbolic form as a three dimensional toric terminal quotient singularity.
\begin{thm}[\cite{Reid},\ (5.2) Theorem]\label{terminal}
A $3$-fold cyclic quotient singularity is terminal if and only if it is of type $\frac{1}{r}(1,a,r-a)$ with $a$ coprime to $r$.
\end{thm}
We note that (ii)-($\beta$) can be expressed as
$$G=\frac{1}{dr}(1,dr-1,ad)$$
where ${\rm gcd}(r,a) = 1,\ r>a$ and $ d>1$. 

We shall compute an irreducible component ${\rm Hilb}^G(\compl{3})$ of $G$-Hilbert scheme $\GHilb{3}$ dominating $\compl{3}/G$ (see Section \ref{G-Hilb}), Fujiki-Oka resolutions (see Section \ref{FO}) and iterated Fujiki-Oka resolutions (see Section \ref{IFO}) for these quotient singularities in each case. These resolutions are toric varieties.

Let $N\cong \inte{n}$ be a free $\inte{}$-module of rank $n\geq 1$, and let $N_{\real{}} := N \otimes_{\inte{}}\real{}\cong \real{n}$. The {\it age} of an element $\nu$ is defined as the sum of all components of $\nu\in N$. Let $\mathcal{H}$ be the affine hyperplane of level $1$,
$$\mathcal{H}:=\{(x_1,x_2,\ldots , x_n)\in N_{\real{n}}\ |\ x_1 +x_2 + \cdots + x_n =1\}.$$ 
Let $\mathfrak{s}$ be the {\it junior simplex} for a rational strongly convex polyhedral cone $\sigma \subset N_{\real{}}$,
$$\mathfrak{s}:=\sigma \cap \mathcal{H}.$$
Moreover, the {\it junior simplex} $\mathfrak{s}$ for a fan $\Delta$ is
$$\mathfrak{s}:=\bigcup_{\sigma\in \Delta}{\left(\sigma\cap \mathcal{H}\right) }.$$
Every $\sigma \subset N_{\real{}}$ has the {\it Hilbert basis} of $\sigma$ with reference to $N$
$$
\mathbf{Hilb}_{N}(\sigma):=\left\{ \nu \in \sigma \cap (N \backslash \{ 0 \}) \left|
\begin{array}{c}
 \text{$\nu$ can not be expresed as}\\  
 \text{the sum of two other vectors}\\ 
 \text{belonging to $\sigma \cap (N\backslash \{ 0 \})$}
\end{array}
\right.\right\}.
$$

Let $\Delta$ be a fan in which all cones are rational strongly convex polyhedral, and let $\Delta(r)$ be the set of all $r$ dimensional cones in $\Delta$. If $\rho \in \Delta(1)$, then there exists a primitive element $P(\rho) \in N \cap \rho$ with $\rho = \mathbb{R}_{\geq0} P(\rho)$. Therefore, we have 
the set of minimal generators of $\sigma \in \Delta$
$$
{\rm Gen}(\sigma):=\{P(\rho) \  | \  \rho \in \Delta(1),\  \rho \prec \sigma \}
$$
where the symbol $\rho \prec \sigma$ means that the face $\rho$ is contained in $\sigma$ as its face. As a natural extension, we also define the set of minimal generators for a fan $\Delta$

$${\rm Gen}(\Delta):=\bigcup_{\sigma \in \Delta} {\rm Gen}(\sigma).$$

\begin{defi}\upshape
The subdivision $\Delta$ of $\sigma$ is called a {\it Hilbert desingularization} of $\sigma$ if $\Delta$ satisfies the following conditions:
\noindent
\begin{itemize}
\item $\Delta$ is smooth,
\item ${\rm Gen}(\Delta) = \mathbf{Hilb}_{N}(\sigma)$.
\end{itemize}
Sometimes, this desingularizetion is simply written as $\HHilb$-{\it desingularization}.
\end{defi}
For the details of this definition , see Section 4 in \cite{Dais}. The previous works on $\HHilb$-desingularizations related to this paper are as follows.
\noindent
\begin{itemize}
\item For two dimensional toric singularity (i.e. cyclic quotient singularity $\compl{2}/G$), there exists a uniquely determined minimal resolution, and this resolution is a $\HHilb$-desingularization.
\item In dimension three, C. Bouvier and G. Sprinberg shows existence of $\HHilb$-desingularization. However, it is not a unique\cite{BS}.
\item They give an example of non-existence of such desingularizations in dimension four\cite{BS}.
\item  A crepant resolution is one of $\HHilb$-desingularizations for toric quotient singularities in any dimension\cite{DHZ}.
\item For three dimensional terminal quotient singularities, Danilov \cite{Danilov} and Reid \cite{Reid} introduce economic resolutions which is obtained by the sequence of weighted blow-ups. It coincides with a $\HHilb$-desingularization.
\end{itemize}
Additionally, the goal of this paper is the following theorem.

\textbf{Theorem \ref{thm1}.}\ 
{\it
For any three dimensional canonical cyclic quotient singularity, there exists a Hilbert iterated Fujiki-Oka resolution.
}\\

\section{On $G$-Hilbert schemes $\Hilb^G(\compl{n})$}\label{G-Hilb}

We recall definition of $G$-graph. Let $S=\mathbb{C}[x_1,\dots,x_n]$ denote the coordinate ring of $\mathbb{C}^n$ and $\mathcal{M}$ be the set of all monomials in $S$ and $1$. Let $\rho_i$ be irreducible representation of $G$. The symbol $X^u$ denotes a monomial $x_1^{u_1}\cdots x_n^{u_n}$ in $\mathcal{M}$ where $u=(u_1, \dots, u_n) \in \mathbb{Z}_{\geq0}^n$.\par
We write ${\rm wt}(X^u)=\rho_i$ if $X^u(g\cdot p)=\rho_i(g)X^u(p)$ holds for any $g \in G$ and $p \in \mathbb{C}^n$. Since any monomial is contained in some $\rho_i$, we can define a map ${\rm wt}: \mathcal{M} \to {\rm Irr}(G)$, where ${\rm Irr}(G)$ is the set of irreducible representation of $G$.
In this paper, we define $G$-graph by using the map "wt" for an ideal in $S$.

\begin{table}[hbtp]
	\begin{tabular}{|l||c|c|} \hline
		     & $SL(n,\mathbb{C})$ & $GL(n,\mathbb{C})$ \\ \hline \hline
		n=2& \multicolumn{2}{|c|}{a minimal resolution} \\
		     & (Ito, Nakamura\cite{IN}) & (Kidoh\cite{Kidoh}, Ishii\cite{Ishii}) \\ \hline
		n=3& a crepant resolutions & Singular in general \\
		     &  (Nakamura\cite{Nakamura}, Bridgeland, King, Reid\cite{BKR})& Example\cite{HIS}: $G=\frac{1}{4}(1,2,3)$ \\ \hline
		n=4& Not crepant in general  & \\
 		     & Example in the terminal case: & Singular in general \\
		     & $G=\frac{1}{2}(1,1,1,1)$ & \\ \hline
	\end{tabular}\\
\begin{center}
Table 1: Properties of ${\rm Hilb}^G(\mathbb{C}^n)$ for $\mathbb{C}^n/G$.
\end{center}
\end{table}
\begin{defi}\label{ggraph}
{\rm Let $I \subset S$ be an ideal, we define a subset $\Gamma(I) \subset \mathcal{M}$ such that $\{X \in \mathcal{M} \mid X \notin I \} $.
A Subset $\Gamma(I)$ is called a} $G$-graph {\rm if the restriction map ${\rm wt}: \Gamma(I) \to {\rm Irr}(G)$ is a bijection.}
\end{defi}

In other words, the condition $\{X \in \mathcal{M} \mid X \notin I \} $ is equivalent to $X \in \Gamma$ and $X$ is divided by $Y \in \mathcal{M}$, then $Y \in \Gamma$. This ideal is called a {\it definition ideal}.

\begin{defi}\label{sigma}
{\rm Let $A_{\Gamma}$ be a set of minimal generators of $I(\Gamma)$. 
We difine the map ${\rm wt}_{\Gamma}:\mathcal{M} \to G\mathchar`-{\rm graph}$ as ${\rm wt}_{\Gamma}(X^u)=\tilde{X}^u$ such that ${\rm wt}(X^u)={\rm wt}(\tilde{X}^u)$. For a $G$-graph $\Gamma$, we define $S(\Gamma)$ to be the subsemigroup of $M$ generated by $\frac{m\cdot x}{{\rm wt}_{\Gamma}(m\cdot x)}$ for all $m \in \mathcal{M}$ and $x \in \Gamma$.
In addition, we define the rational cone
$$
\sigma(\Gamma):=\{w \in N_{\mathbb{R}} \mid w \cdot X^u > w \cdot {\rm wt}_{\Gamma}(X^u) \  for\  all\  X^u \in A_{\Gamma}\},
$$
where $w \cdot X^u$ means standard inner product $w \cdot u$ in $\mathbb{R}^n$.  }
\end{defi}
$\mathrm{Fan}(G)$ denotes the fan defined by all $n$-dimensional closed cone $\sigma(\Gamma)$ and all their faces in $N_{\mathbb{R}}$. The following theorem says that we can calculate an irreducible component ${\rm Hilb}^G(\compl{n})$ of $G$-Hilbert schem $\GHilb{n}$ dominating $\compl{n}/G$ by using $G$-graph.

\begin{thm}\label{Nakamura}{\rm (\cite{Nakamura}, Theorem 2.11)} The following hold.
\begin{itemize}
\item $\mathrm{Fan}(G)$ is a finite fan with its support $\Delta$.
\item The normalization of ${\rm Hilb}^G(\compl{n})$ is isomorphic to the toric variety determined by $\mathrm{Fan}(G)$.
\end{itemize}
\end{thm}

For $G = \frac{1}{r}(1,a,b)$, let $v_1=\frac{1}{r}(1,a,b) \in N$  and $v_i=\frac{1}{r}(i,\bar{ai},\bar{bi})$ where $\bar{x} \equiv x\ (\rm{ mod} \  r)$.\par

The case of $G = \frac{1}{9}(1,4,7)$. The fan ${\rm Fan}(G)$ consists 21 pieces of three dimensional cones. Let $u_1 = e_1 + v_3$, $u_2 = e_2 + v_3$ and $u_3 = e_3 + v_3$. Then these three points are not in Hilbert basis. It follows that ${\rm Gen}({\rm Fan}(G)) \neq \mathbf{Hilb}_{N}(\sigma)$.

\begin{figure}[hbtp]
  \begin{center}
   \includegraphics[width=195pt]{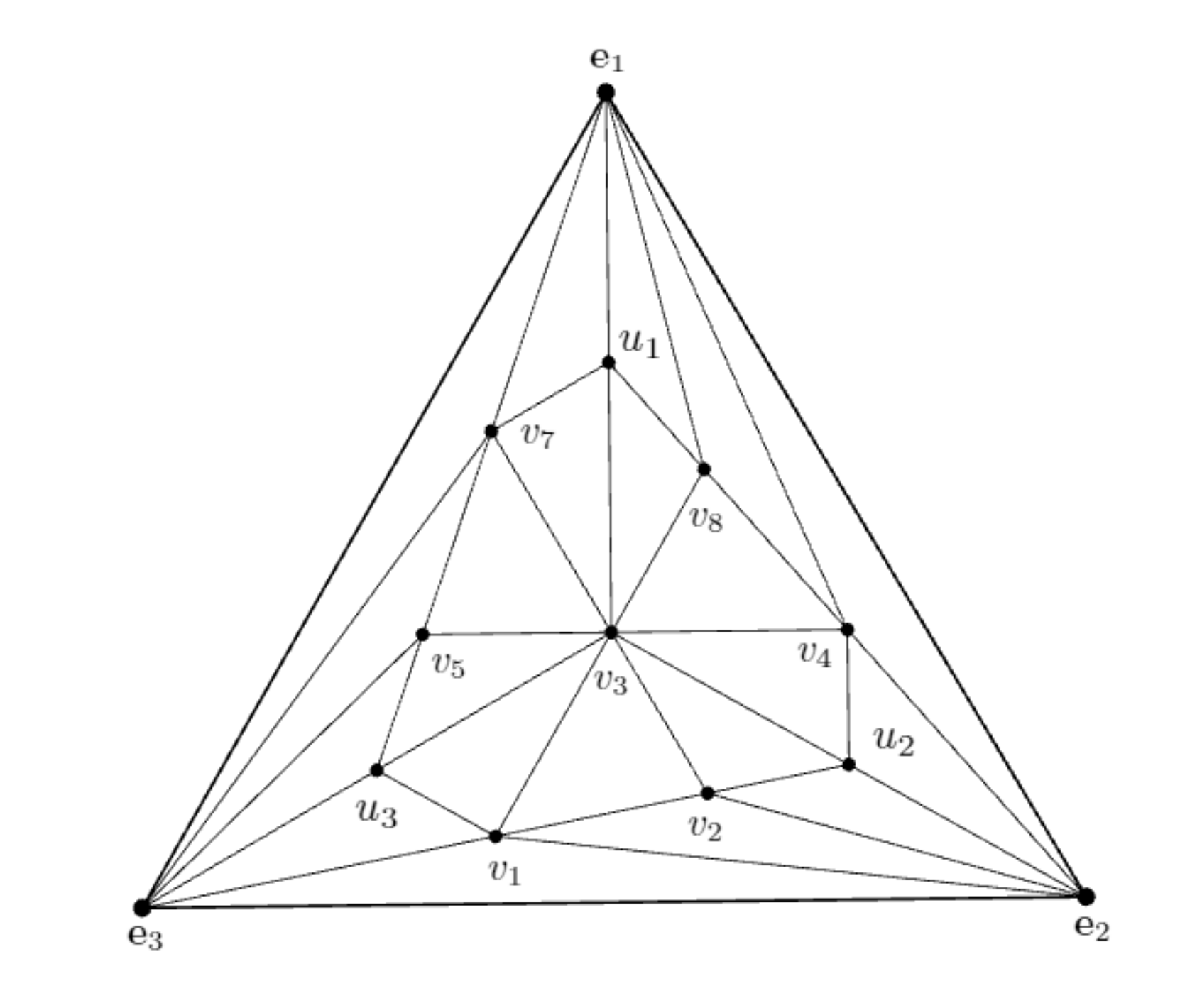}
  \end{center}
  \caption{$\mathfrak{s}$ of $\mathrm{Fan}(G)$ for $\frac{1}{9}(1,4,7)$.}
  \label{147}
\end{figure}

Since $\GHilb{3}$ is reducible in the case of $\frac{1}{14}(1,9,11)$ (see \cite{CMT}), the toric variety ${\rm Hilb}^G(\compl{3})$ is different from $\GHilb{3}$.

\begin{figure}[hbtp]
  \begin{center}
   \includegraphics[width=170pt]{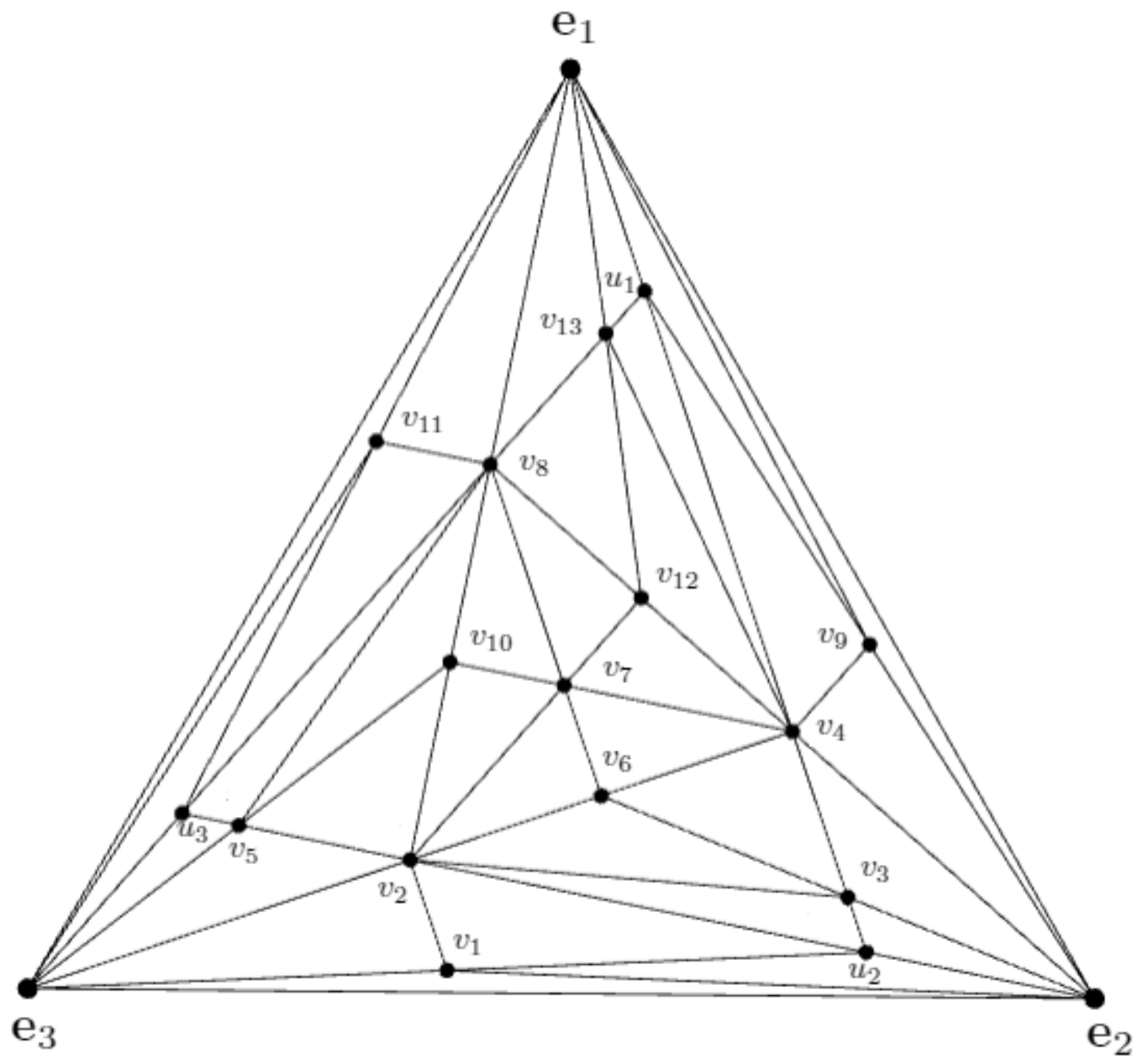}
  \end{center}
  \caption{$\mathfrak{s}$ of $\mathrm{Fan}(G)$ for $\frac{1}{14}(1,9,11)$.}
  \label{147}
\end{figure}
The following six elements in ${\rm Gen}({\rm Fan}(G))$ are not in $ \mathbf{Hilb}_{N}(\sigma)$.
$$v_6 = v_2 + v_4,\ v_{10} = v_2 + v_8,\ v_{12}= v_4 + v_8,$$
$$u_1 = e_1 + v_4,\ u_2 = e_2 + v_2 , u_3 = e_3 + v_8.$$
This tirangulation is not a $\HHilb$-desingularization.

\begin{prop}
Let $G=\frac{1}{4k}(1,2k+1,4k-2)$ with $k > 1$. Then ${\rm Hilb}^G(\compl{3})$ is singular.
\end{prop}

\proofname
For the ideal $I=(x^2, y^2, xy, z^{2k}, xz^k, yz^k)$, the subset $\Gamma(I)$ is
$$
\Gamma(I) =\{1, z, z^2 \dots, z^{2k-1}, x ,xz , xz^{2} \dots, xz^{k-1}, y, yz, yz^2 \dots, yz^{k-1}  \}.
$$
In this case, the map ${\mathrm wt}$ satisfies 
$${\rm wt}(z) = \rho_{4k-2},\ {\rm wt}(xz) = \rho_{4k-1},\ {\rm wt}(yz) = \rho_{2k-3}.$$
Since we have the restriction map as
$$
{\rm Irr}(G) = \{1,\rho_{4k-2}, \rho_{4k-4},\dots, \rho_{2}, \rho_1, \rho_{4k-1}, \dots, \rho_{4k-2(k-2)-1}, \rho_{2k+1}, \rho_{2k-3}, \dots, \rho_3 \},
$$
$\Gamma(I)$ is a $G$-graph.

For the above $\Gamma(I)$, we have the subsemigroup
$$S(\Gamma(I)) = \left\langle \frac{x^2}{z^{2k-1}}, \frac{y^2}{z^{2k-1}}, \frac{xz^k}{y}, \frac{yz^k}{x} \right\rangle ,$$
and the corresponding toric variety $X(N,\sigma(\Gamma(I)))$ satisfies
$$X(N,\sigma(\Gamma(I))) \cong \mathbb{C}[X,Y,Z,W]/(XW-YZ).$$
This is singular.
\qed

\begin{prop}
Let $G=\frac{1}{rd}(1,rd-1,ad)$ with ${\rm gcd}(r,a)=1$. Then $\Hilb^G(\compl{3})$ is singular.
\end{prop}

\proofname
We fix $\Gamma \subset \mathcal{M}$ as 
$$
\Gamma= \{ 1, x, \dots x^{ad-1}, z , y , y^2 , \dots y^{rd-ad-1} \}.
$$
In this case, the subsemigroup $S(\Gamma)$ is
$$
S(\Gamma) = \left\langle \frac{x^{ad}}{z}, \frac{y^{rd-ad}}{z}, \frac{yz}{x^{ad-1}}, \frac{xz}{y^{rd-ad-1}} \right\rangle.
$$
Therefore, the corresponding toric variety satisfies
$$
X(N,\sigma(\Gamma)) \cong \mathbb{C}[X,Y,Z,W]/(XZ-YW),
$$
and this is singular.
\qed

Therefore, we have the following table.

\begin{table}[hbtp]
 \begin{tabular}{|l||c|c|} \hline
   & $\Hilb^G(\compl{3})$  \\ \hline \hline
  ${\rm (i)}$ & a crepant resolution $=$ a $\HHilb$-desingularization \\ \hline
(ii)-($\alpha$) & Singular (O.Kedzierski\cite{Kedzierski}) \\
\hspace{0.5cm} ($\beta$) & Singular  \\  \hline
(iii)& Singular \\ \hline
(iv)& Not $\HHilb$-desingularizations  \\ \hline

  \end{tabular}\\
\begin{center}
Table 2: Hilbert property of $\Hilb^G(\compl{3})$
\end{center}
\end{table}

\section{On Fujiki-Oka resolutions}\label{FO}

The quotient singularity whose corresponding cone in $N_{\mathbb R}$ has at least one smooth facets as an affine toric variety is said to be {\it semi-isolated}. In other words, a semi-isolated quotient singularity is the one containing at least one smooth lines $\compl{n}$ as affine toric subvarieties. The Fujiki-Oka resolution is a canonical resolution for semi-isolated quotient singularities. For the details, see \cite{Fujiki}, \cite{Oka}. This resolution can be constructed algorithmically by the computation of the {\it Ashikaga's continued fraction} for some {\it proper fractions}.

\begin{defi}\upshape\label{profrac}
Let $n$ be an integer greater than or equal to $1$. Let $\mathbf{a}=(a_1,\dots,a_n) \in \inte n$ and $r \in \natu{}$ which satisfies $0\leq  a_i \leq r-1$ for $1\leq i \leq n$. We call the symbol
$$
\frac{\mathbf{a}}{r}=\frac{(a_1,\dots,a_n)}{r}
$$
an {\it $n$-dimensional proper fraction}. Moreover, the proper fraction such that at least one components of $\mathbf{a}$ are $1$ is said to be {\it semi-unimodular}.
\end{defi}

\begin{note}\upshape \label{note}
For the cyclic groups $G$ in the cases (ii)-($\beta$), (iii), and (iv) in the classification of three dimensional canonical cyclic quotient singularities in Section \ref{Intro}, the generator $g$ of $G$ is semi-unimodular if and only if 
\begin{itemize}
\item[(ii)-($\beta$)] $g=\frac{1}{dr}(1,dr-1,ad)$ or $\frac{1}{dr}(dr-1,1,dr-ad)$,
\item[(iii)] $g=\frac{1}{4k}(1,2k+1,4k-2)$ or $\frac{1}{4k}(2k+1,1,4k-2)$,
\item[(iv)] $g=\frac{1}{9}(1,4,7)$, $\frac{1}{9}(7,1,4)$, $\frac{1}{9}(4,7,1)$, $\frac{1}{14}(1,9,11)$, $\frac{1}{14}(11,1,9)$ or $\frac{1}{14}(9,11,1)$
\end{itemize}
where ${\rm gcd}(r,a) = 1,\ r>a,\ d>1$ and $k\geq2$.
\end{note}

The symbol $\overline{\ratio{prop}_n}$ means the set of $n$-dimensional proper fractions and the formal element $\infty$. Similarly, $\overline{\inte n}:= \inte n \cup \{\infty\}$.

Ashikaga's continued fraction consists of a {\it round down polynomial} and a {\it remainder polynomial}, and these polynomials are obtained via {\it round down maps} and {\it remainder maps} for a semi-unimodular proper fraction. For the details, see \cite{Ashikaga}. In this paper, the round down polynomials play important roles especially.

\begin{defi}\upshape
Let $\frac{(1,a_2,\dots,a_n)}{r}$ be a semi-unimodular proper fraction. For $2\leq i \leq n$, the {\it $i$-th remainder map} $R_i:\overline{\ratio{prop}_n} \to \overline{\ratio{prop}_n}$ is defined by
$$
R_i\left(\frac{(1, a_2,\dots,a_n)}{r}\right):=
\left\{ 
\begin{array}{cc}
 \frac{\left(\overline{1}^{a_i}, \overline{a_2}^{a_i},\ \dots,\ \overline{a_{i-1}}^{a_i},\ \overline{-r}^{a_i},\ \overline{a_{i+1}}^{a_i},\ \dots, \overline{a_n}^{a_i}\right)}{a_i}
         & {\rm if} \  a_i \neq 0\\  
 \infty &  {\rm if}\  a_i=0
\end{array}
\right.
$$
and $R_i(\infty)=\infty$ where $\overline{a_j}^{a_i}$ is an integer satisfying $0\leq \overline{a_j}^{a_i} < a_i$ and $\overline{a_j}^{a_i} \equiv a_j$ modulo $a_i$.

\end{defi}

\begin{defi}\upshape \label{DOACF}
Let $\frac{\mathbf{a}}{r}$ be an $n$-dimensional semi-unimodular proper fraction. The {\it remainder polynomial} $\mathcal{R}_*\left(\frac{\mathbf{a}}{r}\right) \in \overline{\ratio{prop}_n}[x_2,\dots,x_n]$ is defined by
  $$
  \mathcal{R}_*\left(\frac{\mathbf{a}}{r}\right):=\frac{\mathbf{a}}{r}+
                                      \sum_{(i_1,i_2,\dots,i_l)\in \mathbf{I}^l,\: l\geq 1 }(R_{i_l}\cdots R_{i_2}R_{i_1})\left(\frac{\mathbf{a}}{r}\right)\cdot x_{i_1}x_{i_2}\cdots x_{i_l}
  $$
  where we exclude terms with coefficients $\infty$ or $\frac{(0,0,\dots,0)}{1}$.

\end{defi}

\begin{ex}\upshape
Let $X(N^{\prime},\sigma)$ have a quotient singularity of $\frac{1}{11}(1,2,8)$-type, i.e., $N^{\prime}=\inte{3}+\inte{}\frac{1}{11}(1,2,8)$ and $\sigma=\real{}_{\geq 0}{e}_1+\real{}_{\geq 0}{e}_2+\real{}_{\geq 0}{
e}_3$. Then, the cone $\sigma$ is semi-unimodular over ${e}_1$, and the Oka center is $c=\frac{1}{11}(1,2,8)$, and the remainder polynomial of the proper fraction $\frac{(1,2,8)}{11}$ is
\begin{eqnarray}
\mathcal{R}_*\left(\frac{(1,2,8)}{11}\right)=
             \frac{1}{11}(1,2,8)&+&\frac{1}{2}(1,1,0)x_2+\frac{1}{8}(1,2,5)x_3 \nonumber \\
                                &+&\frac{1}{2}(1,0,1)x_3x_2+\frac{1}{5}(1,2,2)x_3x_3 \nonumber \\
                                &+&\frac{1}{2}(1,1,0)x_3x_3x_2+\frac{1}{2}(1,0,1)x_3x_3x_3.\nonumber
\end{eqnarray}
This expanding of Ashikaga's continued fraction indicates that the toric variety after the blow-up with the Oka center $\frac{1}{11}(1,2,8)$ has two semi-isolated quotient singularities of $\frac{1}{2}(1,1,0)$-type and $\frac{1}{8}(1,2,5)$-type. For these quotient singularities, the corresponding cones which appear in $\sigma$ after the subdivision by $\frac{1}{11}(1,2,8)\in N^{\prime}$ are $\sigma_2=\real{}_{\geq 0}{e}_1+\real{}_{\geq 0}c+\real{}_{\geq 0}{e}_3$ and $\sigma_3=\real{}_{\geq 0}{e}_1+\real{}_{\geq 0}{e}_2+\real{}_{\geq 0}c$ respectively. $\frac{1}{2}(1,1,0)$ and $\frac{1}{8}(1,2,5)$ are the Oka center of semi-unimodular cones $\sigma_2,\ \sigma_3$ over ${e}_1$ respectively. Therefore, we can take blow-ups with the Oka centers again. The blow-up with Oka centers of $X(N^{\prime},\sigma_3)$ consists a smooth toric variety and quotient singularities of $\frac{1}{2}(1,0,1)$-type and $\frac{1}{5}(1,2,2)$-type respectively. By repeating blow-ups with Oka centers, we have the smooth toric variety (see Fig. \ref{128}). 

\begin{figure}[h]
  \begin{center}
   \includegraphics[width=312pt]{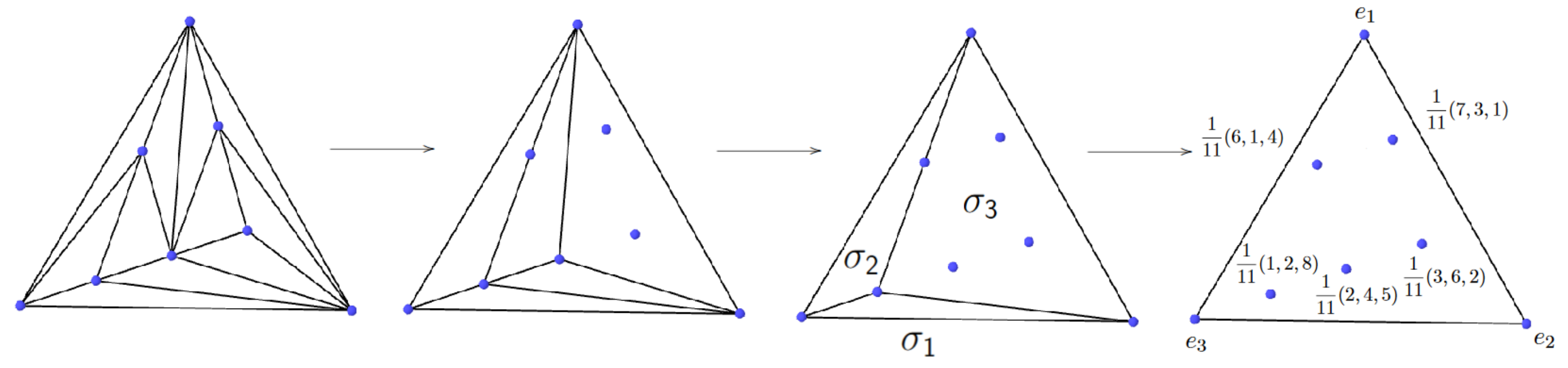}
  \end{center}
  \caption{ $\mathfrak{s}$ of the basic triangulation by Fujiki-Oka resolution}
  \label{128}
\end{figure}

\end{ex}

We proved that all Fujiki-Oka resolutions are crepant for any three dimensional semi-isolated Gorenstein quotient singularity in \cite{SS}. 
\begin{cor}\label{maincor}
For all three dimensional semi-isolated Gorenstein quotient singularities, the Fujiki-Oka resolutions are crepant.
\end{cor}
By Theorem 6.1 in \cite{DHZ}, it follows that arbitrary crepant toric resolution for a Gorenstein abelian quotient singularity is a $\HHilb$-desingularization. Therefore, arbitrary crepant Fujiki-Oka resolution is $\HHilb$-desingularization for any three dimensional semi-isolated Gorenstein quotient singularity. This is the semi-isolated part of (i)  in Section \ref{Intro}.

By Section 5 in \cite{Reid}, the economic resolution for the quotient singularity of $\frac{1}{r}(1,a,r-a)$-type where $\gcd(r,a)=1$ (i.e., the three dimensional toric terminal quotient singularity) is a $\HHilb$-desingularization. For that quotient singularity, the economic resolution coincides with the Fujiki-Oka resolution. Therefore, arbitrary Fujiki-Oka resolution is $\HHilb$-desingularization even in the case of (ii)-($\alpha$). 

For the quotient singularity of $\frac{1}{4k}(1,2k+1,4k-2)$-type where $k\geq2$ (resp. $\frac{1}{dr}(1,dr-1,ad)$-type where ${\rm gcd}(r,a) = 1,\ r>a$ and $ d>1$), the number of the smooth maximal cones appearing in the toric fan of each Fujiki-Oka resolution is greater than or equal to $2|G|$ by the computation of the remainder polynomial of the proper fractions of the both types. We note that this property does not depend on changing semi-unimodular generator of $G$ by Note \ref{note}. Let $E$ and $F$ be the numbers of the edges and faces in $\mathfrak{s}$ of the toric fan corresponding to the Fujiki-Oka resolution respectively. In these case, we have the equation

$$E=\dfrac{1}{2}(3F+s+2)$$
where $s=\gcd(4k,4k+2)=2$ (resp. $s=d$). We note that $s>1$. Let $V$ be the the number of the vertices in $\mathfrak{s}$. By Euler's polyhedron formula, we have
$$V=\dfrac{1}{2}(F+s+4).$$
Since the above inequality $F\geq2|G|$ and $s>1$, we have
$$V\geq|G|+3.$$
On the other hand, the following inequality with respect to the number of the Hilbert basis $ \mathbf{Hilb}_{N}(\sigma)$
\begin{equation}\label{Hilb}
 \#\mathbf{Hilb}_{N}(\sigma) \leq |G| +2
\end{equation}
holds. Therefore, we have
$$V\geq \#\mathbf{Hilb}_{N}(\sigma)+1.$$
This inequality says that the Fujiki-Oka resolution is not $\HHilb$-desingularization in each case of (ii)-($\beta$) and (iii).

For the quotient singularities of $\frac{1}{9}(1,4,7)$-type and $\frac{1}{14}(1,9,11)$-type, the numbers of the smooth maximal cones appearing in the toric fan of the Fujiki-Oka resolution are 23 and 43 respectively by the computation of the remainder polynomial of the proper fractions. The numbers of the smooth maximal cone are invariant under changing semi-unimodular generator of $G$ by Note \ref{note}. By Euler's polyhedron formula and the relation of the number of faces and edges in the triangulation of simplex spanned by $e_1, e_2, e_3$ corresponding to the Fujiki-Oka resolution, the formula
$$V=\dfrac{1}{2}(F+5)$$
holds. Therefore, the number of the vertices $V$ in the simplex is $14$ and $24$ in each case. Clearly, the inequality (\ref{Hilb}) does not stand, and the Fujiki-Oka resolutions in the case (iv) are not $\HHilb$-desingularization. 

\begin{table}[hbtp]
 \begin{tabular}{|l||c|c|} \hline
   & Fujiki-Oka resolution \\ \hline \hline
  $\text{s.i. part of (i)}$ & a crepant resolution $=$ a $\HHilb$-singularization   \\ \hline
(ii)-($\alpha$) &  a $\HHilb$-desingularization (ecnomic) \\
\hspace{0.5cm} ($\beta$) &  Not $\HHilb$-desingularizations \\  \hline
(iii)&  Not $\HHilb$-desingularizations \\ \hline
(iv)&  Not $\HHilb$-desingularizations \\ \hline

  \end{tabular}\\
\begin{center}
Table 3: Hilbert property of Fujiki-Oka Resolutions
\end{center}
\end{table}

\section{On iterated Fujiki-Oka resolutions}\label{IFO}

It is known that there exists an {\it iterated Fujiki-Oka resolution} for Gorenstein abelian quotient singularities in all dimensions by Lemma 4.2 in \cite{SS}. Naturally, the definition of this resolution can be extended to the cyclic quotient singularities. Let $G\subset \glmc{n}$ be a cyclic subgroup and $H$ be a component of a decomposition by cyclic subgroups of $G$. If the singularity $\compl{n}/H$ is semi-isolated, then we have the Fujiki-Oka resolution $(\widetilde{Y_H}, \mathrm{FO}_1)$ and the toric partial resolution $(Y_G,\phi)$ satisfying the following diagram:

\[
\xymatrix@R=6pt{
                              &&&  \compl{n}\ar[dd]^{\pi_H} \\ 
                            && &                 \\ 
                     \widetilde{Y_H} \ar[rrr]^(.5){\mathrm{FO}_1}_{\textrm{ Fujiki-Oka\ Resolution}}\ar[dd]_{\pi_{G/H}}       &&& \compl{n}/H \ar[dd]^{\pi_{G/H}}  \\ 
                          &\circlearrowright &  &  \\ 
                     Y_G=\widetilde{Y_H}\Big/(G/H) \ar[rrr]^{\phi}_{{\rm\qquad  Toric\ Partial\ Resolution}}  &&& \compl{n}/G                  }
\]
where $\pi_{H}$ (resp. $\pi_{G/H}$) is the quotient map by $H$ (resp. $G/H$). If all maximal cones in the toric fan of $Y_G$ are semi-unimodular, then we have a Fujiki-Oka resolutions $(\widetilde{Y_G},\mathrm{FO}_2)$ for the quotient singularities in $Y_G$.
\[
\xymatrix@R=6pt{
\widetilde{Y_G} \ar[rrrr]^{\mathrm{FO}_2}_{\textrm{ Fujiki-Oka\ Resolution}}   &&&&  Y_G }\]

We call the resolution $(\widetilde{Y_G}, \mathrm{FO}_2\circ \phi)$ an {\it iterated Fujiki-Oka resolution} for a cyclic quotient singularity $\compl{n}/G$. Clearly, an alternative Fujiki-Oka resolution can be seen an iterated Fujiki-Oka resolution.

We proved the existence of a crepant iterated Fujiki-Oka resolution for any three dimensional Gorenstein abelian quotient singularity \cite{SS}.

\begin{cor}\label{cor2}
Assume that $G$ is a finite abelian subgroup of $\slmc{3}$.
Then, a crepant iterated Fujiki-Oka resolution exists for $\compl3/G$.
\end{cor}

For the case of (iii) and (iv), we have convenient subgroup to find Hilbert iterated Fujiki-Oka resolutions.

\begin{lem}\label{subg}
Let $G$ be $\frac{1}{4k}(1,2k+1,4k-2)$ with $k > 1$, $\frac{1}{9}(1,4,7)$ or $\frac{1}{14}(1,9,11)$. Then, $G$ contains a cyclic subgroup in $\slmc{3}$.
\end{lem}
\proofname
If $G=\frac{1}{4k}(1,2k+1,4k-2)$ with $k > 1$, then $G$ contains
$$H_1:=\frac{1}{2k}(1,1,2k-2).$$
If $G=\frac{1}{9}(1,4,7)$, then $G$ contains 
$$H_2:=\frac{1}{3}(1,1,1).$$ 
If $G=\frac{1}{14}(1,9,11)$, then $G$ contains
$$H_3:=\frac{1}{7}(1,2,4).$$
\qed

Since the generators $H_1,\ H_2$ and $H_3$ are semi-unimodular, we have Gorenstein cyclic quotient singularities as $\compl{n}/H$ in the above diagram.

\begin{thm}\label{thm1}
For any three dimensional canonical cyclic quotient singularity, there exists a Hilbert iterated Fujiki-Oka resolution.
\end{thm}
\proofname
Since a toric crepant resolution is a $\HHilb$-desingularizations by the discussion in Section \ref{FO}, arbitrary crepant iterated Fujiki-Oka resolution for Gorenstein cyclic quotient singularities is a $\HHilb$-desingularization. This is in the case (i). The Hilbert property for the case of (ii)-($\alpha$) has been shown in Section \ref{FO}. The remaining part is in the case of (iii) and (iv). 

By Lemma \ref{subg}, we have the Fujiki-Oka resolutions $(\widetilde{Y_{H_i}},{\rm FO}_i)$ where $i=1,2,3$. Then, $Y_G$ has a three dimensional toric terminal quotient singularity in each affine piece. By Theorem \ref{terminal}, arbitrary three dimensional toric terminal quotient singularity is $\frac{1}{r}(1,a,r-a)$-type for integers $r$ and $a$ where $0<a<r$ and ${\rm gcd}(r,a)=1$.
We note that the Hilbert basis in each affine piece is not necessary a part of the Hilbert basis globally. However, the age of elements which are not in Hilbert basis globally are greater than or equal to $2$. Therefore, this discussion can be reduct to the case (i) in each affine piece, and we take $\HHilb$-desingularizations as ${\rm FO}_2$.
Then, we have an iterated Fujiki-Oka resolution $(\widetilde{Y_G}, \mathrm{FO}_2\circ \phi)$ which is a $\HHilb$-desingularization.
\qed

By the above theorem, we have the following table.
\begin{table}[hbtp]
 \begin{tabular}{|l||c|c|} \hline
   &iterated Fujiki-Oka resolution \\ \hline \hline
  ${\rm (i)}$ & a crepant resolution $=$ a $\HHilb$-desingularization   \\ \hline
(ii)-($\alpha$) &  a $\HHilb$-desingularization (ecnomic) \\
\hspace{0.5cm} ($\beta$) &   $^\exists$ a $\HHilb$-desingularization \\  \hline
(iii)&  $^\exists$ a $\HHilb$-desingularization \\ \hline
(iv)&  $^\exists$ a $\HHilb$-desingularization \\ \hline

  \end{tabular}\\
\begin{center}
Table 4: Hilbert property of iterated Fujiki-Oka Resolutions
\end{center}
\end{table}



\begin{thebibliography}{}
\bibitem{Ashikaga} T. Ashikaga, {\it Multidimensional continued fractions for cyclic quotient singularities and Dedekind sums}, Kyoto J. Math. (4) {\bf 59} (2019), no.4, 993--1039.
\bibitem{BS} C. Bouvier and G. Sprinberg, {\it Syst\`{e}me g\'{e}n\'{e}rateur minimal, diviseurs essentiels et $G$-d\'{e}singularisations de vari\'{e}t\'{e}s toriques}, Tohoku Math. J. (2) {\bf 47} (1995), no.1, 125--149.
\bibitem{BKR} T. Bridgeland, A. King and M. Reid, {\it The McKay correspondence as an equivalence of derived categories}, J. Amer. Math. Soc. {\bf 14} (2001), no. 3, 535--554. 
\bibitem{CMT} A. Craw, D. Maclagan and R. Thomas, {\it Moduli of McKay quiver representations II: Grobner basis techniques}, J. Algebra {\bf 316} (2007), no. 2, 514--535.
\bibitem{CR} A. Craw and M. Reid, {\it How to calculate $A$-Hilb($\mathbf{C}^3$)}, Geometry of toric varieties, 129--154, Sémin. Congr., {\bf 6}, Soc. Math. France, Paris, 2002.
\bibitem{Dais} D. I. Dais, {\it Resolving $3$-dimensional toric singularities}, Geometry of toric varieties, 155--186, S\'{e}min. Congr., {\bf 6}, Soc. Math. France, Paris, 2002.
\bibitem{DHZ} D. I. Dais, M. Henk, and G. M. Ziegler, {\it On the existence of crepant resolutions of Gorenstein abelian quotient singularities in dimensions $\geq$ 4}, Algebraic and geometric combinatorics, 125--193, Contemp. Math., {\bf 423}, Amer. Math. Soc., Providence, RI, 2006.
\bibitem{Danilov} V. I. Danilov, {\it Birational geometry of three-dimensional toric varieties}, Izv. Akad. Nauk SSSR Ser. Mat. {\bf 46} (1982), no. 5, 971--982, 1135. 
\bibitem{Fujiki} A. Fujiki, {\it On resolutions of cyclic quotient singularities}, Publ. Res. Inst. Math. Sci. {\bf 10} (1974/75), no. 1, 293--328.
\bibitem{HIS} T. Hayashi, Y. Ito and Y. Sekiya, {\it Existence of crepant resolutions}, Higher dimensional algebraic geometry—in honour of Professor Yujiro Kawamata's sixtieth birthday, 185--202, Adv. Stud. Pure Math., {\bf 74}, Math. Soc. Japan, Tokyo, 2017.
\bibitem{II} M. Ishida and N. Iwashita, {\it Canonical cyclic quotient singularities of dimension three}, Complex analytic singularities, 135--151, Adv. Stud. Pure Math., {\bf 8}, North-Holland, Amsterdam, 1987. 
\bibitem{Ishii} A. Ishii, {\it On the McKay correspondence for a finite small subgroup of GL(2,C)}, J. Reine Angew. Math. {\bf 549} (2002), 221--233.
\bibitem{IN} Y. Ito and I. Nakamura, {\it Hilbert schemes and simple singularities}, New trends in algebraic geometry (Warwick, 1996), 151--233, London Math. Soc. Lecture Note Ser., {\bf 264}, Cambridge Univ. Press, Cambridge, 1999.
\bibitem{Jung} S. J. Jung, {\it Terminal Quotient Singularities in Dimension Three via Variation of GIT}, J. Algebra {\bf 468}  (2016), 354--394.
\bibitem{Kedzierski} O. Kedzierski, {\it The $G$-Hilbert scheme for $\frac{1}{r}(1,a,r-a)$}, Glasg. Math. J. {\bf 53} (2011), no. 1, 115--129.
\bibitem{Kidoh} R. Kidoh, {\it Hilbert schemes and cyclic quotient singularities}, Hokkaido Math. J. {\bf 30} (2001), no. 1, 91--103.
\bibitem{Nakamura} I. Nakamura, {\it Hilbert schemes of abelian group orbits},  J. Algebraic Geom. {\bf 10} (2001), no. 4, 757--779.
\bibitem{Oka} M. Oka, {\it On the resolution of hypersurface singularities}, Complex analytic singularities, 405--436, Adv. Stud. Pure Math., {\bf 8}, North-Holland, Amsterdam, 1987. 
\bibitem{Reid} M. Reid, {\it Young person's guide to canonical singularities}, Algebraic geometry, Bowdoin, 1985 (Brunswick, Maine, 1985), 345--414, Proc. Sympos. Pure Math., {\bf 46}, Part 1, Amer. Math. Soc., Providence, RI, 1987. 
\bibitem{SS} K. Sato and Y. Sato, {\it Crepant Property of Fujiki-Oka Resolutions for Gorenstein Abelian Quotient Singularities}, arXiv:2004.03522, preprint (2020).
\end{thebibliography}
\end{document}